\documentclass[preprint,12pt,times]{article}
\usepackage{graphicx}
\usepackage{amsthm}
\usepackage{amssymb}
\usepackage{amsmath}
\usepackage{amscd}
\usepackage{bbm}
\usepackage{float}
\usepackage{booktabs}
\numberwithin{equation}{section}

\newtheorem{theorem}{Theorem}[section]
\newtheorem{lemma}{Lemma}[section]

\usepackage[numbers,sort&compress]{natbib}

\numberwithin{figure}{section}
\numberwithin{table}{section}

\newcommand{\ep}{\mathrm{e}}
\allowdisplaybreaks[2]

\newcommand\btd{\raise 2pt \hbox{$\hat\bigtriangledown$}\hskip 1.5pt}
\newcommand\bt{\raise 2pt \hbox{$\bigtriangledown$}\hskip 1.5pt}

\begin{document}
\date{}
\title{Wave breaking in the unidirectional non-local wave model}
\author{Shaojie Yang\thanks{Corresponding author:
shaojieyang@kust.edu.cn (Shaojie Yang); 19855615281@163.com (Jian Chen)},~~Jian Chen\\~\\
\small ~ Department of Systems Science and Applied Mathematics, \\
 \small~Kunming University of Science and Technology,  \\
 \small ~Kunming, Yunnan 650500, China}

\date{}
\maketitle
\begin{abstract}

In this paper we study wave breaking in the unidirectional non-local wave model describing the motion of a collision-free plasma in a magnetic field.  By analyzing the monotonicity and continuity properties of a  system of the Riccati-type differential inequalities involving the extremal slopes of flows,   we show a new sufficient condition on the initial data to exhibit wave breaking.  Moreover,  the estimates of life span and wave breaking rate are derived. \\


\end{abstract}
\noindent\rule{15.5cm}{0.5pt}

\section{Introduction}\label{sec1}
The motion of a cold plasma in a magnetic field consisting of singly-charged particles can be
described by \cite{r1,r2}
\begin{align}
\label{eq1} n_t+(un)_x=0,\\
\label{eq2} u_t+uu_x+\frac{bb_x}{n}=0,\\
\label{eq3} b-n-\left( \frac{b_x}{n} \right)=0,
\end{align}
where $n=n(t,x), u=u(t,x)$ and $b=b(t,x)$ stand for the ionic density, the ionic velocity and the magnetic field, respectively.
System \eqref{eq1}-\eqref{eq3} can also be applied to describe the motion of collision-free two-fluid model  under the assumptions that the electron inertia, the displacement current and the charge separation are neglected and the Poisson equation \eqref{eq3} is satisfied initially \cite{r3,r4}.

Recently,  Alonso-Or\'{a}n D, Dur\'{a}n A and Granero-Belinch\'{o}n R \cite{r23} derived a unidirectional asymptotic model describing the motion of a collision-free plasma in a magnetic field given in system \eqref{eq1}-\eqref{eq3}.
The method to derive the unidirectional asymptotic model relies on a multi-scale expansion (see \cite{r5,r6,r7,r8,r9,r10,r11}) which reduces the full system \eqref{eq1}-\eqref{eq3} to a cascade of linear equations which can be closed up to some order of precision. More precisely,  let
\begin{align*}
n=1+N,~~~U=u,~~~b=1+B,
\end{align*}
introducing the following formal expansions
\begin{align}\label{eq4}
N=\sum_{l=0}^\infty\varepsilon^{l+1}N^{(l)},~~~B=\sum_{l=0}^\infty\varepsilon^{l+1}B^{(l)},~~~U=\sum_{l=0}^\infty\varepsilon^{l+1}U^{(l)}
\end{align}
then under the  extra assumption $U^{(0)}=N^{(0)}$ in \eqref{eq4} leads to the following bidirectional single non-local wave equation
\begin{align}\label{eq5}
h_{tt}+\mathcal{L}h= (hh_x+[ \mathcal{L},   \mathcal{N}h]h)_x-2(hh_t)_x.
\end{align}
The formal reduction of \eqref{eq5} to the corresponding unidirectional version, cf.\cite{r12}, yields
\begin{align}\label{eq6}
h_t+\frac{3}{2}hh_x=\frac{1}{2}( [ \mathcal{L},   \mathcal{N}h]h+ \mathcal{N}h+h_x  ),
\end{align}
where
\begin{align*}
\mathcal{L}=-\partial_x^2(1-\partial_x^2)^{-1},~~~~\mathcal{N}=\partial_x(1-\partial_x^2)^{-1}
\end{align*}
and $[\mathcal{L}, \cdot   ]\cdot$ denotes the following commutator
\begin{align*}
[\mathcal{L}, f ]g=\mathcal{L}(fg)-f \mathcal{L}g.
\end{align*}

Denote by  $L$ the operator $(1-\partial_x^2)^{-1} $ which acting on functions $f\in L^2(\mathbb{R})$ has the representation
\begin{align*}
Lf(x)=G*f(x)=\int_{\mathbb{R}}G(x-y)f(y)dy,~~~~G(x)=\frac{1}{2}\ep^{-|x|},
\end{align*}
a simple computation, for all $f\in L^2$, we have
\begin{align*}
\partial_x^2Lf(x)=(L-I)f(x),
\end{align*}
where $I$ denotes the identity operator.  Note that $-\mathcal{L}=L-1$, then Eq.\eqref{eq6} can be rewritten as
\begin{align}\label{eq7}
h_t+\frac{3}{2}hh_x=\frac{1}{2}Lh_x+\frac{1}{2}h_x-\frac{1}{2}[L, Lh_x]h.
\end{align}
Eq.\eqref{eq7} possesses the following conservation properties
\begin{align*}
\int_{\mathbb{R}}h(t,x)dx=\int_{\mathbb{R}}h_0(x)dx
\end{align*}
and
\begin{align*}
\int_{\mathbb{R}}h^2(t,x)dx=\int_{\mathbb{R}}h^2_0(x)dx.
\end{align*}

Note that the unidirectional non-local wave model \eqref{eq7} is very similar to the celebrated Fornberg-Whitham equation
\begin{align}\label{eq8}
u_t+\frac{3}{2}uu_x=Lu_x,
\end{align}
which  was derived by Whitham \cite{r13} and Fornberg and Whitham \cite{r14}, as a model for shallow water waves describing breaking waves. Recently,  some properties of FW equation  \eqref{eq8}  such as well-posedness, continuity, traveling waves and blow-up have been studied in \cite{r15,r16,r17,r18,r19,r20,r21,r22}.  The main difference between the unidirectional non-local wave model  and   the FW equation is the emergence of the nonlocal commutator-type term in \eqref{eq7}.

In this paper, we study wave breaking in the unidirectional non-local wave model \eqref{eq7}.
By analyzing wave breaking conditions on the extremal functions
\begin{align*}
  M(t):=\sup\limits_{x\in\mathbb{R}}[h_x(t,x)],
  \end{align*}
and
\begin{align*}
 m(t):=\inf\limits_{x\in\mathbb{R}}[h_x(t,x)],
 \end{align*}
 where both $M(t)$ and $m(t)$ satisfy the Riccati-type differential inequalities   \eqref{ed4} and \eqref{ed5},   we show a new sufficient condition on the initial data to exhibit wave breaking which extending  wave breaking  results in Ref.\cite{r23}.    Our wave breaking condition (see Theorem \ref{th.1}) is that
\begin{align}\label{A}
\inf_{x\in\mathbb{R} }h_{0,x}(x)<\min\left\{-\frac{1}{6}\left(1+\sqrt{1+24C_0}\right),-\frac{1}{12}\left(1+\sqrt{1+24\left(\sup_{x\in\mathbb{R} }h_{0,x}(x)+4C_0\right)}\right)\right\},
\end{align}
where  $C_0$ is some positive constant depends on $\|h_0\|_{L^2}$.  The wave breaking condition in \cite{r23} is that
\begin{align}\label{B}
\inf_{x\in\mathbb{R} }h_{0,x}(x)\leq -C_1,
\end{align}
where $C_1$ is some positive constant depends on $\|h_0\|_{L^2}$ and $\|h_0\|_{L^\infty}$.
Comparing wave breaking condition \eqref{A} and  condition \eqref{B}, it is clear that  the wave breaking condition \eqref{B}  in \cite{r23} follows from our wave breaking condition \eqref{A}. Moreover, for the first time, we get the estimates of life span and wave breaking rate to the unidirectional non-local wave model \eqref{eq7}.

The remainder of this paper is organized as follows. In Section \ref{sec2}, we recall several useful results which are crucial in deriving wave breaking. In Section \ref{sec3}, we present a new wave breaking result and  wave breaking rate.

\section{Preliminaries}\label{sec2}
In this section, we recall several useful results which are crucial in deriving wave breaking.

First of all, we recall the local well-posedness for Eq.\eqref{eq7}.
\begin{lemma}[\cite{r23}]\label{L1}
Let $h_0\in H^s(\mathbb{R})$ with $s>3/2$, then there exists maximal existence time $T>0$ and a unique local solution $h\in C([0, T), H^s(\mathbb{R}))$ of \eqref{eq7}.
\end{lemma}

The blow-up criteria for Eq.\eqref{eq7} can be formulated as follows.
\begin{lemma}[\cite{r23}]
Let $h_0\in H^s(\mathbb{R})$ with $s>3/2$, and let maximal existence time $T>0$ be the lifespan associated to the solution $h$ to  \eqref{eq7} with $h(0,x)=h_0(x)$. Then $h(t,x)$ blows up in finite time $T$ if and only if
\begin{align*}
\int_0^{T}\|h_x(\tau)\|_{BMO}d\tau=\infty.
\end{align*}
\end{lemma}

Next, the following lemma shows that for the maximal time of existence $T> 0$,  the solutions remains bounded.
\begin{lemma}[\cite{r23}]\label{L3}
Let $h_0\in H^s(\mathbb{R})$ with $s>3/2$, and let $T>0$ be the maximal time of existence of the unique solution $h$ of \eqref{eq7} given by Lemma \ref{L1}. Then
\begin{align*}
\sup_{t\in [0,T) }\|h(t)\|_{L^\infty(\mathbb{R})}<\infty.
\end{align*}
\end{lemma}

 Finally, the following classical lemma  is key to construct a  system of the Riccati-type differential inequalities  involving the extremal slopes of flows.
 \begin{lemma}[\cite{r24}]\label{L4}
Let $T>0$ and $u\in C^1([0,T]; H^2)$. Then for every $t\in [0,T)$ there exist at least one point $\xi(t)\in \mathbb{R}$ with
\begin{equation*}
m(t):=\inf_{x\in\mathbb{R}}[h_x(t,x)]=h_x(t,\xi(t)),
\end{equation*}
and the function $m$ is almost everywhere differentiable on $(0,T)$ with
\begin{equation*}
\frac{dm(t)}{dt}=h_{tx}(t,\xi(t)),~\text{a.e.~on}~(0,T).
\end{equation*}
\end{lemma}
 \noindent {\bf Remark.~}The same statement is valid clearly for the supremum function $M(t):=\sup\limits_{x\in\mathbb{R}}[h_x(t,x)]$.

\section{Wave breaking}\label{sec3}
In this section, we present a new wave breaking result and  wave breaking rate. The main result is as follows.
\begin{theorem}\label{th.1}
Let $h_0\in H^s(\mathbb{R})$ with $s>3/2$, and the maximal existence time $T>0$ be the lifespan associated to the solution $h$ to  \eqref{eq7} with $h(0,x)=h_0(x)$. Assume
that the initial data $m_0=m(0)$ and $M_0=M(0)$ satisfy
\begin{equation}\label{b.1}
\inf_{x\in\mathbb{R} }h_{0,x}(x)<\min\left\{-\frac{1}{6}\left(1+\sqrt{1+24C_0}\right),-\frac{1}{12}\left(1+\sqrt{1+24\left(\sup_{x\in\mathbb{R} }h_{0,x}(x)+4C_0\right)}\right)\right\},
\end{equation}
where $C_0=C\|h_0\|_{L^2}^{2}>0$.  Then the corresponding solution $h$ breaks down in the finite time $T$ with
\begin{equation*}
0<T\leq t^*=\frac{2\sqrt3\big(3m_0^2+\frac{1}{2}m_0-2C_0\big)}{3(3m_0^2+m_0-2C_0)\sqrt{3m_0^2-\frac{1}{2}M_0+\frac{1}{2}m_0-2C_0}}.
\end{equation*}
Moreover, for some $x(t)\in\mathbb{R}$, the blow-up rate is
\begin{equation}\label{d.1}
h_x(t,x(t))\sim-\frac{2}{3(T-t)}~~~~as~~t\rightarrow T.
\end{equation}
\end{theorem}

\noindent {\bf Remark.~}The wave-breaking time $t^*$ in Theorem \ref{th.1} could be optimized by more delicate estimate, which can be obtained from replacing
$m^2(t)=-\frac{2}{3}\varphi(t)+\frac{1}{6}(M-m)+\frac{2}{3}C_0\geq-\frac{2}{3}\varphi(t)$ by $m^2(t)\geq-\frac{2}{3}\varphi(t)-\frac{1}{6}m_0+\frac{2}{3}C_0$.
By investigating the proof in Theorem \ref{th.1}, we are able to obtain a better estimate of upper bound on wave-breaking time by
\begin{equation*}
t^*=\frac{\sqrt{3}\big(3m_0^2+\frac{1}{2}m_0-2C_0\big)}{3\sqrt{2C_0-\frac{1}{2}m_0}(3m_0^2+m_0-2C_0)}\ln\frac{\sqrt{6m_0^2-M_0}+\sqrt{4C_0-m_0}}{\sqrt{6m_0^2-M_0}-\sqrt{4C_0-m_0}}.
\end{equation*}

\begin{proof}
Let
\begin{align*}
M(t)=\sup_{x\in\mathbb{R}} h_x(t,x)=h_x(t,\xi_1(t))
\end{align*}
and
\begin{align*}
m(t)=\inf_{x\in\mathbb{R}} h_x(t,x)=h_x(t,\xi_2(t)).
\end{align*}
 Differentiating \eqref{eq7} with respect to $x$ yields
 \begin{align}\label{ed1}
 h_{tx}=-\frac{1}{2}( 3h_x^2+3hh_{xx}+\partial_x [L, Lh_x]h+L_xh_x-h_{xx}  ).
 \end{align}
 Note that the Sobolev embedding $H^{\frac{1}{2}+\epsilon}\hookrightarrow L^\infty $ for $ \epsilon >0$ and
 $L$ is continuous between $H^s(\mathbb{R})$ and $H^{s+2}(\mathbb{R})$ for any $s\in \mathbb{R}$, yields
 \begin{align}\label{ed2}
\notag \big| \partial_x [L, Lh_x]h \big|\leq& \|[L, Lh_x]h\|_{H^{\frac{1}{2}+\epsilon}}\\
\notag \leq& \|Lh_xh\|_{H^{ -\frac{3}{2}+\epsilon    }}+\|(Lh)^2\|_{H^{\frac{3}{2}+\epsilon  }}\\
 \notag\leq& C\|h\|_{L^2}^2\\
 \leq& C\|h_0\|_{L^2}^2:=C_0
 \end{align}
and by observing
\begin{align}\label{ed3}
L_xh_x=G_x*h_x=-\frac{1}{2}\int_{-\infty}^x\ep^{y-x}h_xdy+\frac{1}{2}\int_x^{+\infty}\ep^{x-y}h_xdy.
\end{align}
Combining \eqref{ed2} and \eqref{ed3},  from \eqref{ed1}, it's easy to find that, for a.e. $t\in (0,T)$
\begin{align}\label{ed4}
M'(t)\leq -\frac{3}{2}M^2+\frac{1}{4}(M-m)+C_0,
\end{align}
\begin{align}\label{ed5}
m'(t)\leq -\frac{3}{2}m^2+\frac{1}{4}(M-m)+C_0.
\end{align}
Now, we adapt the method in \cite{r25} in Theorem \ref{th.1} to prove our result. Note that the condition
\eqref{b.1}, together with \eqref{ed5}, we have
\begin{equation*}
m(0)<-\frac{1}{6}\left(1+\sqrt{1+24C_0}\right)~~\text{and}~~m^{\prime}(0)\leq-\frac{3}{2}m^{2}(0)+\frac{1}{4}\big(M(0)-m(0)\big)+C_0<0.
\end{equation*}
We claim that:
\begin{equation*}
m(t)<-\frac{1}{6}\left(1+\sqrt{1+24C_0}\right)~~\text{for~all}~t\in[0,T).
\end{equation*}
Otherwise, there exists $t_0<T$ such that
\begin{equation}\label{b.5}
m(t_0)=-\frac{1}{6}\left(1+\sqrt{1+24C_0}\right)~\text{and}~ m(t)<-\frac{1}{6}\left(1+\sqrt{1+24C_0}\right)~\text{for~all}~ t\in[0,t_0).
\end{equation}
Let
\begin{equation*}
\varphi(t)=-\frac{3}{2}m^{2}+\frac{1}{4}(M-m)+C_0,
\end{equation*}
then from \eqref{ed4} and \eqref{ed5}, a simple computation yields that, for~a.e.~$t\in[0,t_0)$
\begin{align*}
\varphi^{\prime}(t)&=-3m m^{\prime}+\frac{1}{4}(M^{\prime}-m^{\prime}) \\
       &=-\frac{1}{4}(12m+1)m^{\prime}+\frac{1}{4}M^{\prime}\\
      &\leq-\frac{1}{4}(12m+1)\left(-\frac{3}{2}m^{2}+\frac{1}{4}(M-m)+C_0\right)+\frac{1}{4}\left(-\frac{3}{2}M^{2}+\frac{1}{4}(M-m)+C_0\right)\\
      &=~\frac{3}{2}m(3m^2+m-2C_0)-\frac{3}{8}(M+m)^{2}\\
      &<0.
\end{align*}
Taking into account the fact~$\varphi(0)<0$, it follows from the decreasing of $\varphi(t)$ that
\begin{equation*}
\varphi(t)\leq \varphi(0)<0~~\text{for~a.e.}~ t\in[0,t_0).
\end{equation*}
Thus, we have
\begin{equation*}
m^{\prime}(t)\leq \varphi(t)<0~~\text{for~a.e.} ~t\in[0,t_0),
\end{equation*}
which means that
\begin{equation*}
m(t_0)\leq m(0)<-\frac{1}{6}\left(1+\sqrt{1+24C_0}\right).
\end{equation*}
This contradicts with \eqref{b.5}. Consequently, we obtain
\begin{equation}\label{b.6}
m(t)<-\frac{1}{6}\left(1+\sqrt{1+24C_0}\right),~~\varphi(t)\leq \varphi(0)<0,~~~~\text{for~all}~t\in[0,T),
\end{equation}
and
\begin{equation}\label{b.7}
m^{\prime}(t)<0~~~~\text{for~a.e.}~t\in[0,T).
\end{equation}
On the other hand, in view of \eqref{b.6} and \eqref{b.7}, the definition of $\varphi(t)$ implies that
\begin{equation*}
m^2(t)=-\frac{2}{3}\varphi(t)+\frac{1}{6}(M-m)+\frac{2}{3}C_0\geq-\frac{2}{3}\varphi(t)>0,
\end{equation*}
which shows us
\begin{equation}\label{b.8}
-m(t)\geq\sqrt{-\frac{2}{3}\varphi(t)}.
\end{equation}
Conclusively, from the above discussion and \eqref{b.6}, for~a.e.~$t\in[0,T)$
\begin{align*}
\varphi^{\prime}(t)&\leq \frac{3}{2}m(3m^2+m-2C_0)\\
&\leq-\sqrt{6}\frac{3m^2+m-2C_0}{3m^2-\frac{1}{2}M+\frac{1}{2}m-2C_0}(-\varphi)^\frac{3}{2}(t)\\
&\leq-\sqrt{6}\frac{3m^2+m-2C_0}{3m^2+\frac{1}{2}m-2C_0}(-\varphi)^\frac{3}{2}(t).
\end{align*}
A straightforward computation yields
\begin{align*}
\left(\frac{3m^2+m-2C_0}{3m^2+\frac{1}{2}m-2C_0}\right)^{\prime}&=\left(\frac{\frac{1}{2}m}{3m^2+\frac{1}{2}m-2C_0}\right)^{\prime}\\
&=\frac{m^{\prime}(-\frac{3}{2}m^2-C_0)}{\big(3m^2+\frac{1}{2}m-2C_0\big)^2}>0.
\end{align*}
Hence, we have
\begin{equation*}
\varphi^{\prime}(t)\leq-\sqrt{6}\frac{3m_0^2+m_0-2C_0}{3m_0^2+\frac{1}{2}m_0-2C_0}(-\varphi)^\frac{3}{2}(t)~~~~\text{for~a.e.}~t\in[0,T).
\end{equation*}
Denote
\begin{equation*}
\delta=\sqrt{6}\frac{3m_0^2+m_0-2C_0}{3m_0^2+\frac{1}{2}m_0-2C_0}.
\end{equation*}
Apparently, $\delta>0$ and
\begin{equation}\label{b.9}
\left((-\varphi)^{-\frac{1}{2}}(t)\right)^{\prime}=\frac{1}{2}(-\varphi)^{-\frac{3}{2}}(t)\varphi^{\prime}(t)\leq-\frac{\delta}{2}.
\end{equation}
Integrating \eqref{b.9} from $0$ to $t$ yields
\begin{equation}\label{b.10}
-\varphi(t)\geq\left((-\varphi)^{-\frac{1}{2}}(0)-\frac{\delta}{2}t\right)^{-2}.
\end{equation}
Since $-\varphi(t)>0$ for any $t\geq0$, then $\varphi(t)$ and $m(t)$ break down at the finite time
\begin{equation*}
T\leq t^*:=\frac{2}{\delta(-\varphi)^{\frac{1}{2}}(0)}=\frac{2\sqrt3\big(3m_0^2+\frac{1}{2}m_0-2C_0\big)}{3(3m_0^2+m_0-2C_0)\sqrt{3m_0^2-\frac{1}{2}M_0+\frac{1}{2}m_0-2C_0}}.
\end{equation*}
In the case of $T_0=t^*$, in view of \eqref{b.8} and \eqref{b.10}, we have
\begin{equation*}
m(t)\leq-\sqrt{-\frac{2}{3}\varphi(t)}\leq-\frac{1}{\sqrt{\frac{3}{8}}\delta(t^*-t)},
\end{equation*}
which implies that $m(t)\rightarrow-\infty$ as $t\rightarrow T$.

Now, let us give more insight into the wave breaking mechanism  to Eq. \eqref{eq7}. Assume that the corresponding solution $h$ breaks in finite time $T<\infty$ for some $\xi_2(t)\in\mathbb{R}$, then from Lemma \ref{L3},  we have
\begin{equation}\label{c.1}
 m(t)=h_x(t,\xi_2(t))\rightarrow -\infty~~\text{as}~t\nearrow T~~\text{and}~~\sup_{t\in [0,T)}|h(t,\xi_2(t))|<\infty.
\end{equation}
Then it follows from \eqref{ed1} that, for~a.e.~$t\in[0,T)$
\begin{equation}\label{c.2}
m^{\prime}(t)=-\frac{3}{2}m^2(t)-\frac{1}{2}\partial_x [L, Lh_x]h(t,x(t))-\frac{1}{2}L_xh_x(t,x(t)).
\end{equation}
Note that
\begin{align}\label{c.3}
\big| L_xh_x \big| = \big| Lh-h \big| &\leq  \|G*h\|_{L^\infty}+\|h\|_{L^\infty}\nonumber\\
&\leq \|G\|_{L^2}\|h\|_{L^2}+\|h\|_{L^\infty}\nonumber\\
&=\frac{1}{2}\|h\|_{L^2}+\|h\|_{L^\infty}\nonumber\\
&=\frac{1}{2}\|h_0\|_{L^2}+\|h\|_{L^\infty}.
\end{align}
Combing \eqref{c.2}, \eqref{ed2} and \eqref{c.3} leads to
\begin{align}\label{c.4}
-\frac{C_0}{2}-\frac{1}{4}\|h_0\|_{L^2}-\frac{1}{2}\|h\|_{L^\infty}-\frac{3}{2}m^2(t)\leq m^{\prime}(t)\leq-\frac{3}{2}m^2(t)+\frac{1}{2}\|h\|_{L^\infty}+\frac{1}{4}\|h_0\|_{L^2}+\frac{C_0}{2}.
\end{align}

For small $\varepsilon>0$, from \eqref{c.1} we can find $t_0\in(0,T)$  such that $t_0$ is close enough to $T$ in a such
way that
\begin{equation*}
\left\lvert \frac{\frac{1}{2}\|h\|_{L^\infty}+\frac{1}{4}\|h_0\|_{L^2}+\frac{C_0}{2}}{m^2(t)}\right \rvert<\varepsilon
\end{equation*}
on $(t_0,T)$. Therefore, on such interval, one can infer from \eqref{c.4} that, for~a.e.~$t\in(t_0,T)$
\begin{equation*}
-\frac{3}{2}-\varepsilon\leq-\frac{d}{dt}\left(\frac{1}{m}\right)\leq-\frac{3}{2}+\varepsilon.
\end{equation*}
Integrating the above inequalities on $(t,T)$~with~$t\in(t_0,T)$,  we derive the blow-up rate \eqref{d.1} with $x(t)=\xi_2(t)$.  Thus,  the proof of the Theorem \ref{b.1} is complete.
\end{proof}

\section*{Acknowledgment}
This work is supported by Yunnan Fundamental Research Projects (Grant NO. 202101AU070029).

\end{document}